\documentclass{article}[10pt]
\usepackage{theorem,makeidx,amsmath,amsfonts,latexsym}
\theorembodyfont{\upshape}
\newtheorem{theorem}{Theorem}[section]
\newtheorem{prop}{Proposition}[section]

\begin{document}

\title{\large\textbf{{Generating simple random graphs
with prescribed\linebreak degree distribution}}}

\author{Tom Britton \thanks{Department of Mathematics, Stockholm University.}
\and Maria Deijfen \thanks{Department of Mathematics, Stockholm
University, 106 91 Stockholm, Sweden, e-mail: mia@math.su.se, fax:
+46 8 612 6717.} \and Anders Martin-L\"{o}f
\thanks{Department of Mathematics, Stockholm University.}}

\date{May 2006}

\maketitle

\thispagestyle{empty}

\begin{abstract}

\noindent Let $F$ be a probability distribution with support on
the non-negative integers. Four methods for generating a simple
undirected graph with (approximate) degree distribution $F$ are
described and compared. Two methods are based on the so called
configuration model with modifications ensuring a simple graph,
one method is an extension of the classical Erd\H{o}s-R\'{e}nyi
graph where the edge probabilities are random variables, and the
last method starts with a directed random graph which is then
modified to a simple undirected graph. All methods are shown to
give the correct distribution in the limit of large graph size,
but under different assumptions on the degree distribution $F$ and
also using different order of operations.

\vspace{0.5cm}


\vspace{0.5cm}


\end{abstract}

\section{Introduction}

\noindent A graph consists of two sets of objects: a set of
vertices, $\mathcal{V}$, and a set of edges, $\mathcal{E}$. Each
edge corresponds to a pair of vertices and the graph is said to be
\emph{undirected} if these pairs are unordered so that no
directions are associated with the edges. Graphs with no duplicate
edges and no loops -- that is, with at most one edge between each
pair of vertices and with no edges between a vertex and itself --
are called \emph{simple}. Furthermore, a graph is referred to as
random if some kind of randomness is involved in its construction.
In this paper we will consider graphs that are random in that the
edges are generated by random mechanisms. The question at issue
is, given a set of vertices and a probability distribution $F$ on
the non-negative integers, how do we proceed to construct a simple
undirected graph where the degree of a randomly chosen vertex has
distribution $F$?\medskip

\noindent The simplest random graph model is the
Erd\H{o}s-R\'{e}nyi graph, which was introduced in the late 50's
by Paul Erd\H{o}s and Alfr\'{e}d R\'{e}nyi \cite{ER1,ER2}. In fact
they introduced two, closely related, models. Given a set of $n$
vertices and a number $m\leq {n\choose 2}$, the first model,
denoted by $\mathcal{G}_{n,m}$, is defined as the ensemble of
graphs having exactly $m$ edges, each possible graph appearing
with equal probability. The second model, denoted by
$\mathcal{G}_{n,p}$, is obtained by independently adding each one
of the ${n\choose 2}$ possible edges of the graph with some
probability $p$. The distribution of the vertex degree is then
binomial with parameters $n-1$ and $p$ and, if $p$ is scaled by
$1/n$, we get a Poisson distribution in the limit as $n\rightarrow
\infty$. Erd\H{o}s-R\'{e}nyi graphs have been widely studied and
thorough descriptions of the field can be found in
\cite{B2,JLR}.\medskip

\noindent An area that has received a lot of attention within
statistical physics during the last few years is the use of graphs
as models for various types of complex networks; see e.g.\
\cite{DM,N} and the references therein. Examples of network
structures that have been studied are social networks, power
grids, the structure of the internet and various types of
collaboration networks. Typically, this type of networks are very
large, making it impossible to describe them in detail. A natural
approach then is to let the edges representing the connections in
the network be generated by a random procedure designed so that
the resulting graph captures the features of the real-life network
in question as well as possible. Since, as mentioned, the networks
are usually large, it is particularly urgent that the asymptotic
properties of the graph model agree with empirical
observations.\medskip

\noindent An essential characteristic of a graph is the vertex
degree and, in a random graph, this is a random quantity. For
instance, as mentioned above, in the $\mathcal{G}_{n,p}$ model by
Erd\H{o}s and R\'{e}nyi, the degree of a vertex is asymptotically
Poisson distributed. The Erd\H{o}s-R\'{e}nyi graphs have a very
simple and appealing mathematical structure and a lot of work has
been done on the model. However, empirical studies have shown that
the degree distribution in many real-life networks differs
significantly from a Poisson distribution; see e.g. \cite{LEASA}
(human sexual relationships), \cite{FFF} (physical structure of
the internet) and \cite{BA} (movie actor collaboration network).
Complex networks typically have a more heavy-tailed degree
distribution, often specified by some kind of power law, meaning
that the number of vertices with degree $k$ is proportional to
$k^{-\tau}$ for some exponent $\tau>1$. This type of graphs is
often referred to as \emph{scale-free graphs} and there are
important features of such graphs that are missed out if they are
approximated by Erd\H{o}s-R\'{e}nyi graphs; see e.g.\
\cite{BA}.\medskip

\noindent In view of the above, it is important to be able to
generate random graphs with other degree distributions than
Poisson. The aim of this paper is to contribute at this point by
describing a number of algorithms that, given a probability
distribution $F$ (which will later be subject to various
restrictions), produces simple undirected graphs whose vertex
degree is asymptotically distributed according to $F$ (here,
clearly it is required that Supp$(F)\subseteq \mathbb{N}$, where
Supp$(F)$ denotes the support of $F$). To be more precise, given a
set of $n$ vertices and a random mechanism to generate edges
between them, let $p_k^{(n)}$ denote the probability of a randomly
chosen vertex having degree $k$ and write $F=\{p_k;k\geq 0\}$. Our
task is then to design an edge mechanism such that

\begin{itemize}
\item[(i)] $\lim_{n\rightarrow\infty}p_k^{(n)}=p_k$; \item[(ii)]
the resulting graph is simple and undirected.
\end{itemize}

\noindent In all applications mentioned above, the networks are
simple and undirected. Other applications might involve more
complex networks, for instance the link structure of the
world-wide web constitutes a directed graph and bipartite graphs
-- that is, graphs with two types of vertices and edges running
only between unlike types -- are common within sociology. However,
simple undirected networks is indeed an important class in
applications. The present work provides a rigorous treatment of
the asymptotic behavior of the vertex degree in a number of
possible methods for generating such graphs. Some of the
methods/results are new, but an important point of the paper is
also to survey and gather previously known material. To our
knowledge, this type of comparative study with focus on the vertex
degree does not exist previously in the literature.\medskip

\noindent We mention also that there exist algorithms for
generating simple graphs with a given degree \emph{sequence}
$d_1,\ldots ,d_n$; here, $d_i$ is a non-negative integer
specifying the degree of vertex $i$. For example, in McKay and
Wormald \cite{McKW} an algorithm is described which, under certain
conditions, produces a uniformly selected simple undirected graph
with the exact prescribed degree sequence. This is of course even
better than the algorithms studied in the present paper which only
has the correct degree {\it distribution} in the limit as the
number of vertices tends to infinity. However, the price one has
to pay for this is a more complicated algorithm and more
restrictive assumptions about the degree distribution (the McKay
and Wormald criterion corresponds to the degree distribution
having moments of order $4+\epsilon$). In the sequel we are hence
not interested in obtaining a specific degree sequence, but only
in proving that the distribution converges to the intended one in
the limit.\medskip

\noindent The rest of the paper is organized as follows. In
Section 2 we review the well-known configuration model and
describe how it can be used to generate simple graphs with an
arbitrary prescribed degree distribution. Section 3 treats a model
inspired by Chung and Lu \cite{CL1,CL2} that generates simple
graphs with mixed Poisson degree distributions. In Section 4 we
propose a method that is based on the introduction of directed
edges according to a suitably chosen distribution. This method
produces graphs with a degree distribution whose generating
function contains a Poisson factor. Finally, in Section 5 the
methods are discussed and evaluated.

\section{The configuration model}

\noindent The configuration model was defined independently in
\cite{B1} and \cite{W}, both papers being inspired by \cite{BC}.
The model has later been analyzed in \cite{MR1,MR2,NS} for
instance. Given a probability distribution $F$, the model
describes a way to construct an undirected graph on $n$ vertices,
labelled $v_1,\ldots, v_n$, having degree distribution $F$. It is
defined as follows. For each vertex $v_i$, generate a degree $d_i$
independently from a random variable $D$ with distribution $F$ and
attach $d_i$ ``stubs'' to $v_i$. Then join the stubs of all
vertices pairwise completely at random to form edges between the
vertices. To be more precise, first pick two stubs randomly among
all stubs in the graph and join them. Then pick two stubs at
random from the remaining $\sum_1^nd_i-2$ stubs and join them,
etc.\medskip

\noindent A few problems might occur in the construction of a
graph according to this algorithm. The first obvious problem is if
the sum of all degrees, $\sum_id_i$, happens to be an odd number.
In this case there will always be one remaining stub left over in
the pairing algorithm. However, unless $n$ is odd and $F$ is
concentrated to the odd numbers, this problem is easily solved by
either regenerating the degrees until their sum is even or
removing one stub chosen at random. More serious problems arise
when the aim is to generate a {\it simple} undirected graph, that
is, a graph without loops and multiple edges. In the configuration
model, it is clearly possible for a stub of a given vertex $v_i$
with $d_i\geq 2$ to be matched with another one of the stubs of
$v_i$, resulting in an edge from vertex $v_i$ to itself, that is,
a loop. Similarly, two stubs of $v_i$ could by chance be joined
with two stubs of the same other vertex, with the effect that a
multiple edge is created.\medskip

\noindent So what should we do if we insist on the resulting graph
being simple? Two obvious suggestions are (1): to remove loops and
merge multiple edges into single edges in the generated graph to
obtain a simple graph as final product, or (2): to redo the
algorithm until a simple graph occurs by chance. These methods
will be referred to as ``Erased configuration model'' and
``Repeated configuration model'' respectively. Both methods make
the degree distribution somewhat different from the intended one,
but, as we will see, both of them have the right degree
distribution asymptotically under certain moment conditions on the
degree distribution.

\subsection{The erased configuration model}

\noindent Let $F_n=\{p_j^{(n)};j\ge 0\}$ denote the degree
distribution in the erased configuration model with stub
distribution $F=\{p_j;j\ge 0\}$, that is, $p_j^{(n)}$ is the
probability that a randomly selected vertex has degree $j$ in the
erased configuration model on $n$ vertices. Also, write
$N_j^{(n)}$ for the number of vertices having degree $j$ in the
resulting graph.

\begin{theorem}\label{th-config} If $F$ has finite mean, then
\begin{itemize}
\item[\rm{(a)}] $F_n\to F$ as $n\to \infty$;
\item[\rm{(b)}]
$N_j^{(n)}/n \to p_j$ in probability, that is, the empirical
distribution converges in probability to $F$.
\end{itemize}
\end{theorem}

\noindent The proof of this theorem is a bit lengthy, although not
hard. It is found in the Appendix.\medskip

\noindent Theorem \ref{th-config} is not true without the
assumption that the degree distribution $F$ has finite mean, at
least not when the tail decays like $1-F(x)\sim cx^{-\alpha}$
where $\alpha<1$. A heuristic argument for this goes as follows.
For such a decay of the tail, it is known that $\sum_{i=1}^nD_i$
and $\max_{\{1\le i\le n\}}D_i$ are of the same order, see e.g.\
[11, Ch 13.11]. As a consequence, the probability that a stub is
connected to a stub of the maximal node is non-negligible. It
follows that any node with original degree 2 or more has positive
probability (bounded away from 0) to have more than one stub
connected to the maximal node. But whenever this happens, the
degree of the node is decreased in the erased configuration model
and it follows that the new degree distribution will converge to a
distribution stochastically smaller than $F$.\medskip

\noindent In view of the above, there is no hope that the degree
distribution will stay unaffected by the erasing procedure when
the mean is infinite. However, if the degrees are conditioned to
be smaller than $n^{a}$ for some $a\in(0,1)$ it turns out that
Theorem \ref{th-config} remains valid. Indeed, in many
applications it is artificial to include vertices with degree
larger than $n^a$ for some $a\in(0,1)$, and hence it is sometimes
natural with this type of conditioned degrees. Write $F_{n,a}$ for
the degree distribution in the erased configuration model on $n$
vertices where the number of stubs $D^{n,a}_i$ of vertex $v_i$ has
distribution
$$
P(D_i^{n,a}=j)=\frac{P(D_i=j)}{P(D_i\leq n^a)}\quad \textrm{ for
}j=0,\ldots,n^a,
$$
with $D_i\sim F$. Also, let $N^{(n,a)}_j$ be the number of
vertices with degree $j$ in the graph after loops and multiple
edges have been erased.

\begin{theorem}\label{th:config_trunkerad} For any $a\in(0,1)$, we have
\begin{itemize}
\item[\rm{(a)}] $F_{n,a}\to F$ as $n\rightarrow \infty$;
\item[\rm{(b)}] $N_j^{(n,a)}/n \to p_j$ in probability.
\end{itemize}
\end{theorem}

\noindent Except for a few minor modifications, the proof of this
theorem is analogous to the proof of Theorem \ref{th-config}. The
modifications are described in the Appendix.

\subsection{The repeated configuration model}

\noindent The repeated configuration model consists of performing
the configuration model until it produces a simple graph. As
pointed out in \cite{MR1}, it follows from results in \cite{McK}
that the probability of obtaining a simple graph in the
configuration model converges to a strictly positive constant $c$
if the degree distribution has finite second moment. This implies
that a simple graph is then obtained after a geometrically
distributed number of tries. Of course, such a graph might not be
typical for the configuration model. In particular one might
suspect that the number of edges is somewhat smaller than normal,
since there by chance were no multiple edges or loops. However,
below we use the result from \cite{MR1} to show that the resulting
degree distribution converges to the right one provided that it
has finite second moment. In fact, we show the stronger result
that the empirical degree distribution converges to the intended
distribution.\medskip

\noindent Let $F_n =\{p_j^{(n)};j\geq 0\}$ be the degree
distribution of the repeated configuration model on $n$ vertices
with stub distribution $F=\{ p_j; j\ge 0\}$ and write $N_j^{(n)}$
for the number of vertices having degree $j$ in the resulting
graph.

\begin{prop} Assume that $F$ has finite second moment. Then
\begin{itemize}

\item[\rm{(a)}] $F_n\to F$ as $n\to \infty$; \item[\rm{(b)}]
$N_j^{(n)}/n \to p_j$ in probability.
\end{itemize}

\end{prop}

\noindent\textbf{Proof:} We first show part (b). Let $D_1,\ldots,
D_n$ be i.i.d. random variables with distribution $F$ and let
$\tilde p_j^{(n)}=|\{D_i; D_i=j, i=1,\ldots ,n\}|/n$ denote the
empirical distribution of these $n$ variables; here $|\cdot|$
denotes set cardinality. Also, write $S_n$ for the event that the
configuration model on $n$ vertices produces a simple graph. The
empirical distribution of the repeated configuration model is the
same as the distribution of the vector with elements $\tilde
p_j^{(n)}$ conditioned on $S_n$ and we hence have to show that

\begin{equation}\label{eq:rcfm}
P(|\tilde p_j^{(n)}-p_j|>\epsilon\ |\ S_n)\to 0\quad\textrm{as
}n\to\infty\textrm{ for any }\epsilon >0\textrm{ and any }j.
\end{equation}

\noindent Trivially, we have

\begin{eqnarray*}
P\left(|\tilde p_j^{(n)}-p_j|>\epsilon\ |\ S_n\right) & = &
{P\left(|\tilde
p_j^{(n)}-p_j|>\epsilon\ ,\ S_n\right)\over P(S_n)}\\
& \leq & {P\left(|\tilde p_j^{(n)}-p_j|>\epsilon\right)\over
P(S_n)}.
\end{eqnarray*}

\noindent The numerator here converges to 0 by the law of large
numbers and, by the cited result of \cite{MR1}, the assumption
that $F$ has finite second moment implies that $P(S_n)\to c>0$.
Hence (\ref{eq:rcfm}) follows.\medskip

\noindent To show (a), note that, since $0 \le N_j^{(n)}/n\le 1$,
by dominated convergence, the result in (b) implies that
E$[N_j^{(n)}/n] \to p_j$. But
E$[N_j^{(n)}]=\sum_{i=1}^np_j^{(n)}=np_j^{(n)}$, and the desired
result follows.\hfill$\Box$

\section{The generalized random graph}

\noindent In an Erd\H{o}s-R\'{e}nyi graph on $n$ vertices, the
edges are defined by independent Bernoulli random variables
$\{X_{ij}\}_{i<j}$ with $P(X_{ij}=1)=p$, the event $X_{ij}=1$
signifying the presence of an undirected edge between $v_i$ and
$v_j$. By definition, $X_{ji}=X_{ij}$ for $i<j$ and $X_{ii}=0$ for
all $i$. In this section, we consider a model where the
probability $p_{ij}$ of an edge between two vertices $v_i$ and
$v_j$ is allowed to depend on $i$ and $j$. Special cases of this
have been considered in \cite{CL1,CL2,S}. We will show that, if
the probabilities $\{p_{ij}\}$ are picked randomly in a suitable
way, we get a graph with a degree distribution that is easy to
characterize in the limit when $n\rightarrow\infty$; see Theorem
\ref{th:grg}. Also, the degrees of the vertices are
approximatively independent. The model will be referred to as the
generalized random graph.\medskip

\noindent First, we develop the model in more detail. To this end,
let $X=\{X_{ij}\}_{i<j}$ be the array of edge indicators and write
$P(X_{ij}=1)=p_{ij}=1-q_{ij}$. Since the indicators are
independent, the probability density of $X$ is given by
$$
P(X=x)=\prod_{i<j}p_{ij}^{x_{ij}}q_{ij}^{1-x_{ij}}.
$$
Introducing the odds ratios $r_{ij}=p_{ij}/q_{ij}$ and noting that
$p_{ij}=r_{ij}/(1+r_{ij})$ and $q_{ij}=1/(1+r_{ij})$, this can be
written
$$
P(X=x)=\prod_{i<j}(1+r_{ij})^{-1} \prod_{i<j}r_{ij}^{x_{ij}}.
$$
Moreover, if we specialize to the situation where $r_{ij}=u_iu_j$
for some parameters $u=\{u_i\}_{i=1}^n$ with $u_i\geq 0$ and define
$G(u):=\prod_{i<j}(1+u_iu_j)$, we get

\begin{eqnarray}
P_u(X=x) & = & G^{-1}(u)\prod_{i<j}(u_iu_j)^{x_{ij}}\nonumber\\
& = & G^{-1}(u)\prod_iu_i^{d_i(x)},\label{eq:grgslh-fkt}
\end{eqnarray}

\noindent where $d_i(x)$ is the degree of the vertex $v_i$ in the
configuration $x$, that is, $d_i(x):=\sum_jx_{ij}$. This is a
``canonical'' distribution in the sense of statistical mechanics
with sufficient statistics $\{d_i(X)\}$ and from
(\ref{eq:grgslh-fkt}) we see that the conditional distribution of
$X$ given that $\{d_i(X)=d_i\}$ is uniform, that is, all graphs
with a given degree sequence $\{d_i\}$ have the same probability.
This is indeed a nice property of the model, motivating the use of
the parametrization $r_{ij}=u_iu_j$ instead of the one defined by
$p_{ij}=u_iu_j$ used in \cite{CL1,CL2}.\medskip

\noindent To obtain a formula for the joint generating function of
the degree vector $\{d_i(X)\}$, note that, by
(\ref{eq:grgslh-fkt}), we have

\begin{eqnarray*}
\textrm{E}_u\Big[\prod_it_i^{d_i(X)}\Big] & = &
\sum_xP_u(X=x)\prod_it_i^{d_i(x)}\\
& = & G^{-1}(u)\sum_x\prod_i(t_iu_i)^{d_i(x)}.
\end{eqnarray*}

\noindent Since $\sum_xP_u(X=x)=1$, it follows from
(\ref{eq:grgslh-fkt}) that $\sum_x\prod_iu_i^{d_i(x)}=G(u)$, and
hence we get

\begin{eqnarray}
\textrm{E}_u\Big[\prod_it_i^{d_i(X)}\Big] & = &
G^{-1}(u)G(tu)\nonumber\\
& = &
\prod_{i<j}\frac{1+t_iu_it_ju_j}{1+u_iu_j}.\label{eq:grggfkt}
\end{eqnarray}

\noindent Now consider the situation where the parameters
$\{u_i\}$ are suitably scaled random variables, more precisely, we
set $u_i=W_i/\sqrt{n}$, where $\{W_i\}$ are i.i.d. random
variables with finite mean $\mu_{\mbox{{\tiny{$W$}}}}$ and
$W_i\geq 0$. Write $\{D_i\}$ for the degrees of the vertices in
this setting, that is, $D_i=d_i(X)=d_i(X(W))$. The following
theorem specifies the liming distribution of the $D_i$:s.

\begin{theorem}\label{th:grg}
Consider a generalized random graph on $n$ vertices with edge
probabilities defined by $p_{ij}/q_{ij}=W_iW_j/n$, where $\{W_i\}$
are i.i.d. random variables with mean $\mu_{\mbox{{\tiny{$W$}}}}$
and finite moment of order $1+\varepsilon$ for some
$\varepsilon>0$. We have:
\begin{itemize}

\item[\rm{(a)}]The limiting distribution of a degree variable
$D_k$ as $n\rightarrow\infty$ is mixed Poisson with parameter
$W_k\mu_{\mbox{{\tiny{$W$}}}}$. \item[\rm{(b)}] For any $m$, the
variables $D_1,\ldots D_m$ are asymptotically independent.
\end{itemize}
\end{theorem}

\noindent \textbf{Proof:} By taking $t_k=t$, where $0\leq t\leq
1$, and $t_i=1$ for $i\neq k$ in (\ref{eq:grggfkt}), it follows
that
$$
\textrm{E}\big[t^{D_k}\big]=\textrm{E}\bigg[\prod_{i\neq k}
\frac{1+W_iW_kt/n}{1+W_iW_k/n}\bigg].
$$
Using the Taylor expansion $\log(1+x)=x+O(x^2)$, we see that
$$
\prod_i\frac{1+W_iW_kt/n}{1+W_iW_k/n}=\exp\left\{\frac{W_k\sum_i
W_i}{n}(t-1)+R_n\right\},
$$
where $R_n=O(W_k^2\sum_iW_i^2/n^2)$. To estimate $R_n$, note that
$W_i^2\leq \max_{l}\{W_l\}W_i$. The law of large numbers implies
that $\sum_iW_i/n\rightarrow \mu_{\mbox{{\tiny{$W$}}}}$ and, since
the $W_l$:s have finite $1+\varepsilon$-moment, we have that
$\max_{1\leq l\leq n}\{W_l\}/n\rightarrow 0$. It follows that
$R_n$ converges almost surely to 0 as $n\rightarrow\infty$. Hence
$$
\textrm{E}\big[t^{D_k}\big]\rightarrow
\textrm{E}\big[e^{W_k\mu_{\mbox{{\tiny{$W$}}}}(t-1)}\big]
\quad\textrm{as }n\rightarrow\infty,
$$
and part (a) follows. To establish (b), note that, by taking
$t_i=1$ for $i>m$ in (\ref{eq:grggfkt}) and proceeding as in
proving (a), it can be seen that
$$
\textrm{E}\bigg[\prod_{i=1}^mt_i^{D_i}\bigg]\rightarrow\prod_{i=1}^m
\textrm{E}\Big[e^{W_i\mu_{\mbox{{\tiny{$W$}}}}(t_i-1)}\Big]\quad\textrm{as
}n\rightarrow\infty.
$$
Hence the joint generating function of $(D_1,\ldots,D_m)$
asymptotically factorizes into a product of mixed Poisson
generating functions, as desired.\hfill$\Box$\medskip

\noindent Now recall that our task is to generate a simple random
graph with a given degree distribution $F$. According to the above
theorem, if $F$ is mixed Poisson with parameter distribution $Q$
with finite moment of order $1+\varepsilon$, then this can be done
by using the generalized random graph model with i.i.d.\ weights
$\{W_i\}$ distributed according to
$Q/\sqrt{\mu_{\mbox{{\tiny{$Q$}}}}}$, where
$\mu_{\mbox{{\tiny{$Q$}}}}$ denotes the mean of $Q$. As mentioned
in the introduction, the degree distribution in many real-life
networks is heavy-tailed, the probability of a vertex having
degree $k$ being proportional to $k^{-\tau}$ for some exponent
$\tau>1$. It is not hard to see that heavy-tailed mixed Poisson
distributions with this type of power law behavior can be
accomplished by choosing a heavy-tailed parameter distribution
with the desired exponent.\medskip

\noindent In this context it is clearly of interest to know to
what extent Theorem 3.1 is still true if the distribution of
$\{W_i\}$ has a heavy tail such that $E[W_i]=\infty$. Indeed,
according to the theorem below, if we assume that the tail of the
distribution varies regularly in the sense that $P(W_i>w)\sim
cw^{-\alpha}$ as $w\to\infty$, for some $\alpha\in(0,1)$ and some
constant $c>0$, then the distribution of the $D_k$:s is still
mixed Poisson, but with a different scaling and different mixing
distribution. Note that a power law distribution with exponent
$\tau\in(1,2)$ satisfies the tail condition with $\alpha=\tau-1$.

\begin{theorem}\label{th:grg_inf_mean}
Suppose that $\{W_i\}$ are i.i.d.\ with $P(W_i>w)\sim
cw^{-\alpha}$, for some $\alpha\in(0,1)$ and $c>0$, and consider
the generalized random graph with
$p_{ij}/q_{ij}=W_iW_j/n^{1/\alpha}$. Then:

\begin{itemize}
\item[\rm{(a)}]The limiting distribution of a degree variable
$D_k$ is mixed Poisson with parameter $\gamma W_k^{\alpha}$ where
$\gamma =c\int_0^\infty (1+x)^{-2}x^{-\alpha}dx$.

\item[\rm{(b)}]For any $m$, the variables $D_1,\ldots ,D_m$ are
asymptotically independent.
\end{itemize}
\end{theorem}

\noindent \textbf{Remark} As pointed out above, the tail behavior
of a mixed Poisson distribution with a power law parameter is
determined by the parameter. Hence, for a mixed Poisson variable
$Y$ with parameter $\gamma W_k^{\alpha}$, we have that

\begin{eqnarray*}
P(Y\geq y) & \approx & P(W_k^{\alpha}\geq y)\\
& = & P(W_k\geq y^{1/\alpha})\\
& \sim & cy^{-1}
\end{eqnarray*}

\noindent that is, the distribution is a power law with exponent
$\tau=2$.\medskip

\noindent\textbf{Proof of Theorem \ref{th:grg_inf_mean}:} As in
Theorem \ref{th:grg} we have
$$
E\left[t^{D_k}\right]=E\left[ \prod_{i\ne
k}\frac{1+W_iW_kt/n^{1/\alpha}}{1+W_iW_k/n^{1/\alpha}}\right].
$$
Let us first fix the value of $W_k$, say $W_k=w$. Then
\begin{equation}
\Phi_n(w):= E[t^{D_k}| W_k=w]=E\left[\prod_{i\ne k}\varphi\left(
\frac{W_i}{n^{1/\alpha}}\right) \right],\label{phi}
\end{equation}
with $\varphi(x):=(1+xwt)/(1+xw)$. Write $V(x)$ for the
distribution function of $W_i$. Since all $\{W_i\}$ are
independent, we have

\begin{eqnarray*}
\Phi_n(w) & = & \left(\int_0^\infty
\varphi\left(\frac{x}{n^{1/\alpha}}\right)V(dx)\right)^{n-1}\\
& = & \left( 1+\int_0^\infty
(\varphi(x)-1)V(n^{1/\alpha}dx)\right)^{n-1}\\
& = & \left( 1+\int_0^\infty
\varphi'(x)\big(1-V(n^{1/\alpha}x)\big)dx\right)^{n-1},
\end{eqnarray*}

\noindent where the last equality follows from partial
integration. In order to see that the last integral is $O(1/n)$,
note that
$$
n\int_0^\infty \varphi'(x)\big(1-V(n^{1/\alpha}x)\big)dx=
\int_0^\infty\frac{\varphi'(x)}{x^\alpha}(n^{1/\alpha}x)^\alpha
(1-V(n^{1/\alpha}x))dx.
$$
By the assumption, $y^\alpha(1-V(y))$ is a bounded function which
converges to $c>0$ as $y\to\infty$, and hence, by bounded
convergence,

\begin{eqnarray*}
\lim_{n\rightarrow\infty}\int_0^\infty\frac{\varphi'(x)}{x^\alpha}(n^{1/\alpha}x)^\alpha
(1-V(n^{1/\alpha}x))dx & = & c\int_0^\infty
\frac{\varphi'(x)}{x^\alpha}dx\\
& = & (t-1)w^\alpha\gamma,
\end{eqnarray*}

\noindent where $\gamma:=c\int_0^\infty(1+x)^{-2}x^{-\alpha}dx$.
It follows that $\lim_{n\to\infty}\Phi_n(w)= e^{(t-1)\gamma
w^\alpha}$, which is recognized as the generating function of a
Poisson distribution with mean $\gamma w^\alpha$. Integrating over
$w=W_k$, we see that the limit distribution of $D_k$ is mixed
Poisson with mean $\gamma W_k^{\alpha}$ and generating function
$$
E\left[ t^ {D_k}\right]=E\left[e^{(t-1)W_k^\alpha}\right]=
\int_0^\infty e^{(t-1)w^\alpha}V(dx).
$$
The proof of (b) is analogous to the proof of Theorem \ref{th:grg}
and is therefore omitted.\hfill$\Box$\medskip

\noindent We finish this section by showing that the empirical
degree distribution in the generalized random graph converges to
the asymptotic mixed Poisson degree distribution in the graph. To
this end, when $\{W_i\}$ have finite mean
$\mu_{\mbox{{\tiny{$W$}}}}$, write $N_k^{(n)}$ for the number of
vertices having degree $k$ in the generalized random graph with
edge probabilities defined by $p_{ij}/q_{ij}=W_iW_j/n$, and let
$F=\{p_k;k\geq 0\}$ be a mixed Poisson distribution with parameter
$W\mu_{\mbox{{\tiny{$W$}}}}$. Similarly, when $P(W_i>w)\sim
cw^{-\alpha}$ for some $\alpha\in(0,1)$, the number of vertices
having degree $k$ for edge probabilities defined by
$p_{ij}/q_{ij}=W_iW_j/n^{1/\alpha}$ is denoted by
$N_k^{(n,\alpha)}$, and we write $F^\alpha=\{p_k^\alpha;k\geq 0\}$
for a mixed Poisson distribution with parameter $\gamma
W^{\alpha}$, where $\gamma$ is defined in Theorem
\ref{th:grg_inf_mean}.

\begin{prop}\label{prop:grgef} As $n\rightarrow\infty$ in the
generalized random graph, we have:

\begin{itemize}
\item[(a)] If $\{W_i\}$ have finite moment of order
$1+\varepsilon$, then $N_k^{(n)}/n\rightarrow p_k$ in probability
for all $k$.

\item[(b)] If $P(W_i>w)\sim cw^{-\alpha}$ for some
$\alpha\in(0,1)$, then $N_k^{(n,\alpha)}/n\rightarrow
p_k^{\alpha}$ in probability for all $k$.
\end{itemize}

\end{prop}

\noindent \textbf{Proof:} The proofs of part (a) and part (b) are
analogous and we give here the proof of (a). Write $P^{(n)}$ for
the probability law of the generalized random graph on $n$
vertices and let $\mathbf{1}_{\{\cdot\}}$ denote the indicator
function. Clearly $N_k^{(n)}=\sum_{i=1}^n\mathbf{1}_{\{D_i=k\}}$
and hence, using symmetry, it follows that

\begin{eqnarray*}
\textrm{E}\big[N_k^{(n)}\big] & = & \sum_{i=1}^nP^{(n)}(D_i=k)\\
& = & nP^{(n)}(D_1=k).
\end{eqnarray*}

\noindent By Theorem \ref{th:grg} (a), we have
$P^{(n)}(D_1=k)\rightarrow p_k$ as $n\rightarrow\infty$, meaning
that E$[N_k^{(n)}/n]\rightarrow p_k$. The desired result is now
obtained from Chebyshev's inequality if we can show that
Var$(N_k^{(n)}/n)\rightarrow 0$. To do this, note that

\begin{eqnarray*}
\textrm{E}\Big[\big(N_k^{(n)}\big)^2\Big] & = &
\textrm{E}\Big[\sum_i\mathbf{1}_{\{D_i=k\}}^2+\sum_{i\neq
j}\mathbf{1}_{\{D_i=k\}} \mathbf{1}_{\{D_j=k\}}\Big]\\
& = & nP^{(n)}(D_1=k)+n(n-1)P^{(n)}(D_1=k,D_2=k),
\end{eqnarray*}

\noindent where the last equality follows from symmetry. By
Theorem \ref{th:grg} (b), the variables $D_1$ and $D_2$ are
asymptotically independent and hence we have
$P^{(n)}(D_1=k,D_2=k)\rightarrow p_k^2$ as $n\rightarrow \infty$.
Using the formula Var$(X)=\textrm{E}[X^2]-\textrm{E}[X]^2$, it
follows that Var$(N_k^{(n)}/n)\rightarrow 0$ as
$n\rightarrow\infty$ and we are done. \hfill$\Box$\medskip

\section{Directed graphs with removed directions}

\noindent In this section we propose a construction method where
directed edges are introduced according to some distribution $G$.
The directions of the edges are then disregarded and multiple
edges are fused together, producing a graph whose asymptotic
degree distribution $F$ is the convolution of the distribution $G$
and a Poisson distribution with parameter
$\mu_{\mbox{{\tiny{$G$}}}}$, where $\mu_{\mbox{{\tiny{$G$}}}}$
denotes the mean of $G$; see Proposition
\ref{th:poconv}.\medskip

\noindent To describe the construction, write
$\mathcal{V}=\{v_1,\ldots,v_n\}$ for the vertex set. Let $G$ be a
probability distribution with finite mean
$\mu_{\mbox{{\tiny{$G$}}}}$ and Supp$(G)\subset \mathbb{N}$, and
write $\{g_k\}$ for the probabilities associated with $G$. Also,
define $G_n$ via the probabilities
$$
g_k^{(n)} := \left\{
\begin{array}{ll}
                      g_k & \mbox{for }k=0,1,\ldots,n-2;\\
                      \sum_{k=n-1}^\infty g_k & \mbox{for }k=n-1;\\
                      0 & \mbox{for }k\geq n,
                    \end{array}
            \right.
$$
that is, $G_n$ is a truncated version of $G$ with support on
$\{0,1,\ldots,n-1\}$. The graph is constructed as follows:

\begin{itemize}

\item[1.] Associate independently to each vertex $v_i$ a random
variable $Y_i$ with distribution $G_n$, and add to the graph $Y_i$
directed edges pointing out from $v_i$. The vertices to be hit by
the edges starting at $v_i$ are chosen randomly without
replacement from $\mathcal{V}\setminus \{v_i\}$, independently for
all vertices. This defines a directed random graph
$\mathcal{G}_{dir}(n,G)=\{\mathcal{V},\mathcal{E}_{dir}\}$.

\item[2.] To obtain a simple undirected graph
$\mathcal{G}(n,G)=\{\mathcal{V},\mathcal{E}\}$, the directions of
the edges are disregarded and multiple edges are fused together,
that is, an undirected edge between the vertices $v_i$ and $v_j$
is included in $\mathcal{E}$ as soon as at least one of the
directed edges $(v_i,v_j)$ and $(v_j,v_i)$ is present in
$\mathcal{E}_{dir}$.
\end{itemize}

\noindent Let $D_i$ denote the degree of the vertex $v_i$ in the
resulting undirected graph. To find an expression for $D_i$, write
$\mathcal{V}_i^{out}$ for the set of vertices that are hit by
edges pointing out from $v_i$ in $\mathcal{E}_{dir}$ and write
$\mathcal{V}_i^{in}$ for the set of vertices that sends outgoing
edges to $v_i$ in $\mathcal{E}_{dir}$. Define
$Z_i=|\mathcal{V}_i^{in}\cap\neg \mathcal{V}_i^{out}|$ (here,
$\neg$ denotes set complement) so that, in words, $Z_i$ indicates
the number of edges in $\mathcal{E}_{dir}$ pointing at $v_i$ and
starting at vertices that are not hit by outgoing edges from
$v_i$. Some thought reveals that
$$
D_i=Y_i+Z_i.
$$
Clearly all variables $\{Z_i\}$ have the same distribution, which
we denote by $H_{(n,G)}$. Also remember that $Y_i\sim G_n$ for all
$i$. Hence the degree variables $\{D_i\}$ are identically
distributed and we write $F_{(n,G)}$ for their distribution. The
following theorem is the aforementioned characterization of the
asymptotic degree distribution as the convolution of $G$ and a
Poisson distribution with the same mean as $G$.

\begin{theorem}\label{th:poconv} As $n\rightarrow\infty$, we have
\begin{itemize}

\item[\rm{(a)}] $G_n\rightarrow G$; \item[\rm{(b)}]
$H_{(n,G)}\rightarrow {\rm{Po}}(\mu_{\mbox{{\tiny{$G$}}}})$, where
${\rm{Po}}(\mu_{\mbox{{\tiny{$G$}}}})$ is a Poisson distribution
with mean $\mu_{\mbox{{\tiny{$G$}}}}$; \item[\rm{(c)}]
$F_{(n,G)}\rightarrow G*{\rm{Po}}(\mu_{\mbox{{\tiny{$G$}}}})$.
\end{itemize}

\end{theorem}

\noindent \textbf{Proof:} The claim in (a) is immediate from the
definition of $G_n$. To prove (b), fix a vertex $v_k$ and, for
$i\neq k$, let $X_{ik}$ be a 0-1 variable, indicating whether
there is a directed edge from $v_i$ to $v_k$ in
$\mathcal{E}_{dir}$ or not. Since there are $Y_i$ outgoing edges
from $v_i$ and the vertices to be hit by these edges are chosen
randomly without replacement from the $n-1$ vertices in
$\mathcal{V}\backslash \{v_i\}$, we have $ P(X_{ik}=1)=Y_i/(n-1)$.
Also, we have that

$$
Z_k=\sum_{i\not = k ;v_i\in\neg\mathcal{V}_k^{out}}X_{ik}.
$$

\noindent Hence, by conditioning on $\{Y_i\}$, we obtain
$$
\textrm{E}\left[t^{Z_k}\right]= \textrm{E}\left[\prod_{i\not=k
;v_i\in\neg\mathcal{V}_k^{out}}\left(1+\frac{Y_i}{n-1}(t-1)\right)\right],
$$
and, since $\{Y_i\}$ are i.i.d.\ with distribution $G_n$ and
$\left|\neg\mathcal{V}_k^{out}\backslash \{v_k\}\right|=n-1-Y_k$,
it follows that

\begin{eqnarray*}
\textrm{E}\left[t^{Z_k}\right] & = &
\left(1+\frac{\mu_{\mbox{{\tiny{$G_n$}}}}} {n-1}(t-1)\right)^{n-1}
\textrm{E}\left[\left(\frac{n-1}{n-1+\mu_{\mbox{{\tiny{$G_n$}}}}(t-1)}
\right)^{Y_k}\right].
\end{eqnarray*}

\noindent Here $\mu_{\mbox{{\tiny{$G_n$}}}}$ denotes the mean of
the distribution $G_n$. As $n\rightarrow\infty$, the first factor
on the right hand side converges to $e^{\mu_G(t-1)}$, which is
recognized as the moment generating function of a Poisson variable
with parameter $\mu_{\mbox{{\tiny{$G$}}}}$, and the left hand side
converges to 1. Hence part (b) is established.\medskip

\noindent To prove (c), note that
$\textrm{E}\left[t^{Y_k+Z_k}\right]=\textrm{E}\left[t^{Y_k}
\textrm{E}\left[t^{Z_k}|Y_k\right]\right]$. The inner expectation
is calculated as above and we get
$$
\textrm{E}\left[t^{Y_k+Z_k}\right]=
\left(1+\frac{\mu_{G_n}}{n-1}(t-1)\right)^{n-1}
\textrm{E}\left[\left(\frac{t(n-1)}{n-1+\mu_{G_n}(t-1)}\right)^{Y_k}\right].
$$
As $n\rightarrow\infty$, the second factor converges to the
generating function of $G$ and, as pointed out above, the limit of
the first factor is the moment generating function of a Poisson
distribution with mean $\mu_G$.\hfill$\Box$\medskip

\noindent Having proved that the asymptotic distribution of the
vertex degree in $\mathcal{G}(n,G)$ is $G*$Po$(\mu_G)$, the
obvious question is which distributions can arise in this way. We
will not give a full answer to this question -- it is presumably
difficult -- but rather give a few examples of distributions that
can indeed be obtained and also specify how the distribution $G$
of the number of outgoing arrows at the vertices should be chosen
in these cases.\medskip

\noindent \textbf{Power law distribution}\smallskip

\noindent First note that, if $G$ is a power law distribution with
exponent $\tau$, then $G*{\rm{Po}}(\mu_{\mbox{{\tiny{$G$}}}})$
will be so as well. Hence, if we are not interested in the exact
form of the resulting degree distribution, but only in that its
tail decays as a certain power law, then the
$\mathcal{G}(n,G)$-model is clearly applicable.

\noindent \textbf{Poisson distribution}\nopagebreak\smallskip

\noindent The simplest case when the exact form of the resulting
degree distribution is important is the Poisson distribution.
Clearly, a Poisson distribution $F$ with parameter
$\mu_{\mbox{{\tiny{$F$}}}}$ is accomplished by choosing $G$ to be
a Poisson distribution with parameter
$\mu_{\mbox{{\tiny{$F$}}}}/2$.\medskip

\noindent \textbf{Mixed Poisson distribution}\smallskip

\noindent  A mixed Poisson distribution can be obtained as a
limiting degree distribution in $\mathcal{G}(n,G)$ given a certain
condition (\ref{eq:mixcond}) on the law of the parameter. To see
this, assume that $F$ is mixed Poisson with parameter law $Q$ with
finite mean $\mu_{\mbox{{\tiny{$Q$}}}}$. The moment generating
function of $F$ then equals
$$
\psi_{\mbox{{\tiny{$F$}}}}(t)=\int_0^\infty e^{x(t-1)}dQ(x).
$$
Since it should hold that $F=G*$Po$(\mu_{\mbox{{\tiny{$G$}}}})$,
we have
$\psi_{\mbox{{\tiny{$G$}}}}(t)=\psi_{\mbox{{\tiny{$F$}}}}(t)e^{-\mu_G(t-1)}$
and $\mu_{\mbox{{\tiny{$G$}}}}=\mu_{\mbox{{\tiny{$Q$}}}}/2$, and
hence
$$
\psi_{\mbox{{\tiny{$G$}}}}(t)=\int_0^\infty
e^{(x-\mu_{\mbox{{\tiny{$Q$}}}}/2)(t-1)}dQ(x).
$$
To ensure that this is the generating function of a probability
distribution, let $\xi_Q=\inf\textrm{Supp}(Q)$ and assume that

\begin{equation}\label{eq:mixcond}
\xi_Q-\mu_{\mbox{{\tiny{$Q$}}}}/2>0.
\end{equation}

\noindent Then
$$
\psi_{\mbox{{\tiny{$G$}}}}(t)=\int_0^\infty
e^{y(t-1)}d\widetilde{Q}(y),
$$
where $\widetilde{Q}(y)=Q(y+\mu_{\mbox{{\tiny{$Q$}}}}/2)$, that
is, $\widetilde{Q}$ is the distribution $Q$ translated
$\mu_{\mbox{{\tiny{$Q$}}}}/2$ units to the left. This means that
$G$ is mixed Poisson with parameter distribution $\widetilde{Q}$.
Hence, a mixed Poisson distribution with parameter distribution
$Q$ that satisfies $\xi_Q-\mu_{\mbox{{\tiny{$Q$}}}}/2>0$, is
obtained as a limiting degree distribution in $\mathcal{G}(n,G)$
by choosing $G$ to be mixed Poisson with parameter distribution
$\widetilde{Q}$.\medskip

\noindent \textbf{Compound Poisson distribution}\smallskip

\noindent Let $F$ be compound Poisson with Poisson parameter
$\lambda$ and discrete summand distribution $R$ with finite mean
$\mu_{\mbox{{\tiny{$R$}}}}$ and generating function
$\psi_{\mbox{{\tiny{$R$}}}}(t)$ (recall that a compound Poisson
distribution is the law of a sum of a Poisson number of i.i.d.\
random variables). Then $\mu_{\mbox{{\tiny{$F$}}}}=\lambda
\mu_{\mbox{{\tiny{$R$}}}}$ and $
\psi_{\mbox{{\tiny{$F$}}}}(t)=e^{\lambda(\psi_R(t)-1)}$, and, if
we want the limiting degree distribution in $\mathcal{G}(n,G)$ to
be $F$, then we must have

\begin{eqnarray}
\psi_{\mbox{{\tiny{$G$}}}}(t) & = & \psi_{\mbox{{\tiny{$F$}}}}(t)e^{-\mu_G(t-1)}\nonumber\\
& = & e^{\lambda(\psi_R(t)-1)-\mu_G(t-1)}\label{eq:compgen}.
\end{eqnarray}

\noindent Here, since
$\mu_{\mbox{{\tiny{$G$}}}}=\lambda\mu_{\mbox{{\tiny{$R$}}}}/2$ and
$\psi_{\mbox{{\tiny{$R$}}}}(t)=\sum_0^\infty r_kt^k$, where
$\{r_k\}$ denotes the probabilities associated with $R$, the
exponent in (\ref{eq:compgen}) can be written as
$$
\lambda\left(\left(r_0+\frac{\mu_{\mbox{{\tiny{$R$}}}}}{2}\right)+
\left(r_1-\frac{\mu_{\mbox{{\tiny{$R$}}}}}{2}\right)t+\sum_{k=2}^\infty
r_kt^k-1\right).
$$
Assume that

\begin{equation}\label{eq:compcond}
r_1>\mu_{\mbox{{\tiny{$R$}}}}/2,
\end{equation}

\noindent and introduce a new distribution $R'$ by defining
$$
r'_k = \left\{
\begin{array}{ll}
                     r_0+\mu_{\mbox{{\tiny{$R$}}}}/2  & \mbox{for }k=0;\\
                     r_1-\mu_{\mbox{{\tiny{$R$}}}}/2  & \mbox{for }k=1;\\
                     r_k & \mbox{for }k\geq 2,
                    \end{array}
            \right.
$$
that is, $R'$ is obtained by transferring the mass
$\mu_{\mbox{{\tiny{$R$}}}}/2$ from the point 1 to the point 0 in
the distribution $R$ (note that a consequence of
(\ref{eq:compcond}) is that $R$ must be chosen so that
$\mu_{\mbox{{\tiny{$R$}}}}<2$). With $R'$ defined in this way we
have
$$
\psi_{\mbox{{\tiny{$G$}}}}(t)=e^{\lambda(\psi_{R'}(t)-1)},
$$
that is, $G$ is a compound Poisson distribution with Poisson
parameter $\lambda$ and summand distribution $R'$. Hence a
compound Poisson distribution $F$ with summand distribution $R$
that satisfies (\ref{eq:compcond}) is obtained as limiting degree
distribution in $\mathcal{G}(n,G)$ by choosing $G$ to be a
compound Poisson distribution with summands distributed according
to $R'$.\medskip

\section{Concluding comments}

\noindent In the present paper, four different ways of generating
simple undirected graphs with a prescribed degree distribution are
described. The methods are referred to as the erased configuration
model, the repeated configuration model, the generalized random
graph and the directed graph with removed directions (DGRD)
respectively. None of the methods is able to produce a graph that
has the desired distribution exactly -- that is, in a finite
graph, a randomly selected vertex will not have exactly the
correct degree distribution -- but under certain regularity
assumptions, it is shown that all four methods give the right
distribution in the limit as the number of vertices $n$ tends to
infinity.\medskip

\noindent Let us summarize the assumptions on the degree
distribution for the different methods: In order for the repeated
configuration model to produce a simple graph in stochastically
bounded time as $n\to \infty$, the second moment of the degree
distribution has to be finite and for the generalized random graph
model with edge probabilities scaled by $n$ to be applicable,
finite moment of order $1+\varepsilon$ for some $\varepsilon>0$ is
required. For the other two methods, finite mean is sufficient. If
the degrees are conditioned on being smaller than $n^a$ for some
$a\in(0,1)$, the erased configuration model can handle
distributions with infinite mean as well, and, with a different
scaling of the edge probability, also the generalized random graph
can be applied to infinite mean distributions.\medskip

\noindent As for the class of achievable distributions, the erased
configuration model and the repeated configuration model are both
able to generate graphs with any limiting distribution. The
generalized random graph model can only produce mixed Poisson
distributions and the DGRD-model gives distributions that can be
expressed as the convolution of a discrete distribution with
finite mean and a Poisson distribution with the same mean -- a
class containing certain types of mixed Poisson and compound
Poisson distributions for instance. However, if only tail
properties of the desired distribution are specified, both the
generalized random graph model and the DGRD-model can do the
job.\medskip

\noindent Concerning the number of operations needed to produce
the graph, it is easily seen to be of order $n$ for all methods
except the generalized random graph, which requires $O(n^2)$
operations. Note however that an approximation of the generalized
random graph model that uses only $O(n)$ operations can be
obtained by replacing the conditional degree distribution of a
given vertex, conditional on the $W_i$:s, with a Poisson
distribution with the same mean. Among the methods using $O(n)$
operations the erased configuration model and the DGRD-model
require less operations: The repeated configuration model
generates a graph a geometrically distributed number of times
whereas the other two methods only generates a graph once and then
erases a few edges.\medskip

\noindent A perhaps more subjective opinion is that the
generalized random graph model is probabilistically more tractable
than the other methods. Its construction is straightforward,
containing less dependence structures, implying that it is easier
to show property results for this model. Also, it was easily seen
from the construction of the general random graph that the
obtained graph is uniform in the sense that all graphs with a
given degree sequence have the same probability.\medskip

\noindent Finally we mention that, apart from the degree
distribution, there is of course a number of other properties of a
graph of which it would also be desirable to have control over,
for instance clustering and path length. Rigorous analyzes of
algorithms incorporating these aspects is to a large extent
lacking and it is a future task to further investigate and develop
methods for generating random graphs specifying such
properties.\bigskip

\noindent \textbf{Acknowledgement} We thank Remco van der Hofstad
for giving the idea of how to strengthen Theorem \ref{th-config}
to cover also the empirical degree distribution, and for
suggesting Theorem \ref{th:config_trunkerad}. Maria Deijfen also
thanks van der Hofstad for interesting and enlightening
discussions on the topic of random graphs.

\section*{Appendix}

\noindent \textbf{Proof of Theorem \ref{th-config}:} First note
that, since $0 \le N_j^{(n)}/n\le 1$, by dominated convergence,
the result in (b) implies that E$[N_j^{(n)}/n] \to p_j$. But
E$[N_j^{(n)}]=\sum_{i=1}^np_j^{(n)}=np_j^{(n)}$, and hence (a)
follows from (b).\medskip

\noindent To show (b), write $\widetilde{N}^{(n)}_j$ for the
number of vertices that has degree $j$ before edges are erased to
make the graph simple. By the law of large numbers
$\widetilde{N}^{(n)}_j/n\rightarrow p_j$ as $n\rightarrow\infty$
and hence we are done if we can show that
$\big(\widetilde{N}^{(n)}_j- N^{(n)}_j\big)/n\rightarrow 0$ in
probability as $n\rightarrow \infty$. Let $M^{(n)}$ be the number
of vertices where at least one stub is removed in the erasing
procedure and note that $\widetilde{N}^{(n)}_j- N^{(n)}_j\leq
M^{(n)}$. Hence, by Markovs inequality, it suffices to show that
E$[M^{(n)}]/n\rightarrow 0$. To do this, for $i=1,\ldots,n$, let
$E_i$ be the number of stubs attached to $v_i$ that are rubbed out
in the erasing procedure and define
$$
M^{(n)}_i = \left\{ \begin{array}{ll}
                      1 & \mbox{if $E_i\geq 1$;}\\
                      0 & \mbox{if $E_i=0$}.
                    \end{array}
            \right.
$$
Also, let $D_i$ be the degree of vertex $v_i$ before loops and
multiple edges have been erased and write $P^{(n)}$ for the
probability law of the erased configuration model on $n$ vertices
(that is, averaged out also over the original degrees $\{D_i\}$).
Since $M^{(n)}=\sum_{i=1}^nM^{(n)}_i$ and $\{M^{(n)}_i\}$ are
equally distributed, we have

\begin{eqnarray*}
\frac{1}{n}\,{\rm{E}}\left[M^{(n)}\right] & = &
\frac{1}{n}\sum_{i=1}^n
{\rm{E}}\left[M^{(n)}_i\right]\\
& = & {\rm{E}}\left[M^{(n)}_1\right]\\
& = & P^{(n)}(E_1\geq 1).
\end{eqnarray*}

\noindent The work lies in proving that $P^{(n)}(E_1\geq
1)\rightarrow 0$ as $n\rightarrow\infty$, or equivalently, that
$P^{(n)}(E_1=0)\rightarrow 1$. Clearly this follows if we show
that

\begin{equation}\label{eq:cfmpt1}
P^{(n)}(E_1=0|D_1=j)\to 1 \textrm{ for all }j\textrm{ as
}n\to\infty.
\end{equation}

\noindent To prove (\ref{eq:cfmpt1}), write $A_j$ for the event
that a given stub belonging to a vertex with $j$ stubs in total
avoids being removed in the erasing procedure. Below we show that

\begin{equation}\label{eq:A_j}
P^{(n)}(A_j)\rightarrow 1\quad\textrm{for all $j$ as }n\rightarrow
\infty.
\end{equation}

\noindent This establishes (\ref{eq:cfmpt1}): After having merged
one of the $j$ stubs of the vertex $v_1$ to a stub belonging to
some other vertex, saving it from being erased, the probability
that a fixed one of the other $j-1$ stubs are erased equals
$P^{(n)}(A_{j-1})$, since the fact that one stub from the other
vertices is no longer available for merging is asymptotically
negligible. This can then be repeated until there are no remaining
stubs of $v_1$ and (\ref{eq:cfmpt1}) follows by noting that
$P^{(n)}(E_1=0|D_1=0)=1$.\medskip

\noindent For the proof to be complete it remains to prove
(\ref{eq:A_j}). To do this, first remember that a stub can be
erased for two reasons: because it forms a loop and because it is
part of a multiple edge. For the sake of completeness we also
include the case when a randomly selected stub is removed if the
total number of stubs is odd. Now, consider a fixed stub belonging
to a vertex $v$ with $j$ stubs in all. Write $A_j^{loop}$ and
$A_j^{mult}$ for the events that the stub is not part of a loop
and a multiple edge respectively and let $A_j^{odd}$ be the event
that the stub is not removed as the randomly selected ``odd''
stub. To estimate the probabilities of these events, we condition
on that the total number of stubs equals $m$ and write $P^{(n)}_m$
for the corresponding conditioned probability measure. If $m$ is
odd, the probability that the stub is removed as the ``odd'' stub
is $1/m$ and, if $m$ is even, the probability is 0. Hence

\begin{equation}\label{eq:A_jodd}
P^{(n)}_m\left(A_j^{odd}\right)\geq 1-\frac{1}{m}
\end{equation}

\noindent For the stub to form a loop it has to be joined to one
of the other $j-1$ stubs of the vertex $v$ and, since clearly the
stub is matched to each one of the other $m-1$ stubs in the graph
with the same probability $1/(m-1)$, this happens with probability
$(j-1)/(m-1)$, that is,

\begin{equation}\label{eq:A_jloop}
P^{(n)}_m\left(A_j^{loop}\right)=1-\frac{j-1}{m-1}.
\end{equation}

\noindent To compute the probability that the stub is not part of
a multiple edge, assume that it does not make up a loop and
condition on the degree $k$ of the vertex $v'$ of the stub to
which it is joined.  Also, number the remaining $j-1$ stubs of the
vertex $v$ in some arbitrary way from 1 to $j-1$ and let $B_j$ be
the event that there are no loops among these stubs. Trivially,
$$
P^{(n)}\left(A_j^{mult}\right)\geq P^{(n)}\left(A_j^{mult}\cap
B_j\right).
$$
For the event $A_j^{mult}\cap B_j$ to happen, none of the
remaining $j-1$ stubs of the vertex $v$ can connect to another
stub of $v$ or to a stub originating from $v'$. Considering the
stubs $1,\dots,j-1$ one in a turn, we see that the probability
that stub number 1 avoids being connected to a stub of $v$ or $v'$
is $(m-j-k)/(m-3)$ (the denominator comes from that two stubs are
already used in the fixed connection between $v$ and $v'$ and the
stub cannot join to itself). Then, given that stub number 1 is
connected to some other vertex, the probability that stub number 2
is so as well is $(m-j-k-1)/(m-5)$, and so on. Furthermore, it is
only possible for $A_j^{mult}\cap B_j$ to happen if $j-1\leq
m-k-j$, since otherwise loops among the $j-1$ stubs on vertex $v$
or multiple edges between $v$ and $v'$ cannot be avoided. Hence,
writing $P^{(n)}_{m,k}$ for the probability measure conditioned on
$m$ and $k$, we have
$$
P^{(n)}_{m,k}\left(A_j^{mult}\cap
B_j\right)=\left\{\begin{array}{ll}

                     {m-k-j\over m-3}{m-k-j-1\over m-5}\cdots {m-k-2j+2\over
                      m-(2j-1)} & \mbox{if }j-1\leq m-k-j;\\
                      0 & \mbox{otherwise,}
                    \end{array}
            \right.
$$
implying that

\begin{equation}\label{eq:A_jmult}
P^{(n)}_{m,k}\left(A_j^{mult}\right)\geq \left({m-k-2j+2\over
m-3}\right)^{j-1}_+,
\end{equation}

\noindent where $r_+=\max\{r,0\}$. Combining (\ref{eq:A_jodd}),
(\ref{eq:A_jloop}) and (\ref{eq:A_jmult}) and using Boole's
inequality, it follows that

\begin{eqnarray*}
P^{(n)}_{m,k}(A_j) & = & P^{(n)}_{m,k}\left(A_j^{odd}\cap
A_j^{loop} \cap A_j^{mult}
\right)\\
& \geq & \left({m-k-2j\over
m}\right)^{j-1}_+-\frac{1}{m}-\frac{j-1}{m-1}.
\end{eqnarray*}

\noindent Removing the conditioning on $m$ and $k$ and denoting
the corresponding random variable $L_n$ and $K_n$ respectively, we
get

\begin{equation}\label{eq:A_jbound}
P^{(n)}\left(A_j\right)\geq \textrm{E}\left[\left(
{\frac{L_n-K_n-2j}{L_n}}\right)^{j-1}_+\right]-\textrm{E}\left[\frac{1}{L_n}\right]
-\textrm{E}\left[\frac{j-1}{L_n-1}\right].
\end{equation}

\noindent To complete the proof, we use dominated convergence to
show that the right hand side of (\ref{eq:A_jbound}) converges to
1 as $n\rightarrow \infty$, establishing (\ref{eq:A_j}). Recall
that $L_n$ is the total number of stubs in the configuration and
$K_n$ is the number of stubs connected to the vertex of a randomly
selected stub. The total number of stubs is a sum of $n$ i.i.d.\
random variables $\{D_l\}$ with distribution $F$ and mean
$\mu_{\mbox{{\tiny{$F$}}}}$, which is finite by assumption. Hence,
the law of large numbers implies that
$L_n/n\rightarrow\mu_{\mbox{{\tiny{$F$}}}}$ almost surely as
$n\rightarrow\infty$. The conditional distribution of $K_n$ given
$\{N_i^{(n)}\}$ is specified by
$P(K_n=i)=iN^{(n)}_i/\sum_rrN_r^{(n)}$, where $N^{(n)}_i$ is the
number of $D_l$:s that equal $i$. Since $N^{(n)}_i/n\rightarrow
p_i$ as $n\rightarrow\infty$ and
$\mu_{\mbox{{\tiny{$F$}}}}<\infty$, it follows that $K_n$
converges in distribution to a proper random variable $K$ with
distribution
$P(K=i)=ip_i/\sum_rrp_r=ip_i/\mu_{\mbox{{\tiny{$F$}}}}$ and hence
$K_n/n\rightarrow 0$ almost surely as $n\rightarrow\infty$.
Combining these two observations, we get that
$$
\left({L_n-K_n-2j\over L_n}\right)^{j-1}_+\rightarrow
1\quad\textrm{a.s. as }n\rightarrow\infty.
$$
Furthermore, since $j$ and $K_n$ are both strictly positive, we
have
$$
0\le \left({L_n-K_n-2j\over L_n}\right)^{j-1}_+ \le 1.
$$
By dominated convergence, it follows that the first term on the
right hand side of (\ref{eq:A_jbound}) converges to 1, and it is
easily seen that the other two terms converge to 0. Hence the
proof is complete.\hfill$\Box$\medskip

\noindent \textbf{Proof of Theorem \ref{th:config_trunkerad}:} Fix
$a\in(0,1)$ and let the degrees be distributed according to
$F_{n,a}$. Clearly, if $F$ has finite mean, the claim follows
immediately from Theorem \ref{th-config}, so assume that $F$ has
infinite mean. The proof of Theorem \ref{th-config} then remains
valid all the way to the estimate (\ref{eq:A_jbound}). We need to
see that the right hand side tends to 1 for all $j$ as
$n\rightarrow \infty$. To this end, let $\tilde{\mu}$ be the
expectation of $F$ conditioned on being smaller than, say, the
smallest possible value in the support plus 10. Then clearly
$P(L_n>\tilde{\mu}n)\rightarrow 1$ as $n\rightarrow\infty$, and,
trivially,
$$
\textrm{E}\left[\left({L_n-K_n-2j\over
L_n}\right)^{j-1}_+\right]\geq
\textrm{E}\left[\left({L_n-K_n-2j\over
L_n}\right)^{j-1}_+\bigg|L_n>\tilde{\mu}n\right]P(L_n>\tilde{\mu}n).
$$
On the event $\{L_n>\tilde{\mu}n\}$, we have $K_n/L_n\leq
n^a/\tilde{\mu}n\rightarrow 0$. Hence, by bounded convergence, the
right hand side above tends to 1 as $n\rightarrow\infty$ and,
since the other two terms in the estimate (\ref{eq:A_jbound})
clearly converge to 0, we are done.\hfill$\Box$

\end{document}